\documentclass[hidelinks,onefignum,onetabnum]{siamonline220329}

\usepackage{algorithmic} 

\usepackage{bbm}

\Crefname{ALC@unique}{Line}{Lines}

\usepackage{lipsum}
\usepackage{amsfonts}
\usepackage{graphicx}
\usepackage{epstopdf}
\usepackage{algorithmic}
\ifpdf
  \DeclareGraphicsExtensions{.eps,.pdf,.png,.jpg}
\else
  \DeclareGraphicsExtensions{.eps}
\fi

\usepackage{enumitem}
\setlist[enumerate]{leftmargin=.5in}
\setlist[itemize]{leftmargin=.5in}

\newsiamremark{remark}{Remark}
\newsiamremark{hypothesis}{Hypothesis}
\crefname{hypothesis}{Hypothesis}{Hypotheses}
\newsiamthm{claim}{Claim}

\headers{A Relaxed Wasserstein Distance Formulation for RMMs}{Keyu Chen, Zetian Wang, and Yunxin Zhang}

\title{A Relaxed Wasserstein Distance Formulation for Mixtures of Radially Contoured Distributions\thanks{Submitted to the editors DATE.
\funding{This research is supported by National Key R\&D Program of China (2024YFA1012401), the Science and Technology Commission of Shanghai Municipality (23JC1400501), and Natural Science Foundation of China (12241103).}}}

\author{Keyu Chen\thanks{School of Mathematical Sciences, Fudan University, Shanghai 200433, China. 
  (\email{chenky21@m.fudan.edu.cn}, \email{wangzt24@m.fudan.edu.cn}).}
\and Zetian Wang\footnotemark[2]
\and Yunxin Zhang\footnotemark[2]\thanks{Shanghai Key Laboratory for Contemporary Applied Mathematics, Laboratory of Mathematics for Nonlinear Science, Fudan University, Shanghai 200433, China. 
	(\email{xyz@fudan.edu.cn}).}}

\usepackage{amsopn}

\newcommand{\dx}{\mathrm{d}x}
\newcommand{\dy}{\mathrm{d}y}

\newcommand{\dr}{\mathrm{d}r}
\newcommand{\norm}[1]{\left\|{#1}\right\|}

\newcommand{\R}{{\mathbb R}}

\ifpdf
\hypersetup{
  pdftitle={A Relaxed Wasserstein Distance Formulation for Mixtures of Radially Contoured Distributions},
  pdfauthor={Keyu Chen, Zetian Wang, and Yunxin Zhang}
}
\fi

\begin{document}

\maketitle

\begin{abstract}
Recently, a Wasserstein-type distance for Gaussian mixture models has been proposed. However, that framework can only be generalized to identifiable mixtures of general elliptically contoured distributions whose components come from the same family and satisfy marginal consistency. In this paper, we propose a simple relaxed Wasserstein distance for identifiable mixtures of radially contoured distributions whose components can come from different families. We show some properties of this distance and that its definition does not require marginal consistency. We apply this distance in color transfer tasks and compare its performance with the Wasserstein-type distance for Gaussian mixture models in an experiment. The error of our method is more stable and the color distribution of our output image is more desirable.
\end{abstract}



\section{Introduction}

In recent years, optimal transport (OT) has attracted a lot of attention in mathematics \cite{EvansL.C1999DEMf}, engineering \cite{6502714}, data science \cite{korotin2023neural} and imaging science \cite{BauerMartin2015Ddmb, KarlssonJohan2017Gsif} . This work was originally considered by Monge \cite{monge1781memoire} and has significantly developed after Kantorovich's work \cite{KantorovichL.V.2006Otto}. The main purpose of optimal transport is to measure the distance of two probability distributions, known as Wasserstein distance, and to find an optimal transport map between the two distributions. In practice, optimal transport is formulated with a linear programming problem. As an aspect of applied mathematics, the computation method plays an important part in OT-related work. The basic solvers of OT are the classic linear programming solvers such as simplex method and inner point method \cite{BoydStephenP2004Co}. They have a complexity of $O(n^3)$, and thus the computation is expensive in large-scale problems. Recently, some fast algorithms have been proposed, enabling OT to deal with such problems, such as Sinkhorn algorithm \cite{2013Sinkhorn}, Greenkhorn algorithm \cite{AltschulerJason2017Ntaa}, and others \cite{GuoWenshuo2020FAfC, DvurechenskyPavel2018COTC, LinTianyi2019Oeot}. 

Another way to reduce the computation is trying to reduce the amount of data. A natural way is that we can sample part of data, compute the Wasserstein distance several times and take the average value. Alternatively, one may use the expectation of Wasserstein distance between empirical distributions to approximate the true Wasserstein distance. This method was first researched in \cite{DudleyR.M.1969TSoM}, which proved that the expectation of the Wasserstein distance between empirical distributions converge in expectation to the true Wasserstein distance and gave a convergence rate. More careful discussion about the convergence behavior can be found in \cite{BolleyFrancois2007Qcif, BoissardEmmanuel2011Sbfc, WEEDJONATHAN2019Saaf}.

However, in practice, we usually not only need the Wasserstein distance but also need to find a map that transports one distribution to another. For instance, while dealing with color transfer task, what we actually need is only the  transfer map. The sample-based method can not handle this task. Another way to reduce complexity is to add some structure to the data, which is easier to operate than discrete distributions. The Gaussian mixture models (GMMs) comes first. It is a popular probability model in statistics and machine learning \cite{HastieTrevor2009EoSL}. Thanks to the EM algorithm \cite{DempsterA.P.1977MLfI}, the parameters of GMMs are easily inferred. In image processing, it is used in a large body of works to represent the patch distributions of images, such as image restoration \cite{2011From, GuoshenYu2012SIPW, TeodoroAfonsoM.A.M.2015Sida} and texture synthesis \cite{GalerneBruno2017SOTi}. Recently, a Wasserstein-type distance in the space of GMMs has been proposed in \cite{DelonJulie2020Awdi}. Based on the explicit formulation of the Wasserstein distance and optimal transport map for Gaussian distributions, they give two ways to define a transport map. Their method is shown useful in color transfer task and texture synthesis task.

The work based on GMMs relies on that GMMs can approximate distributions well. If we view it as a function approximation problem, it is approximating a function with Gaussian basis weighted by normalized positive coefficients. However, Gaussian density function is elliptically contoured, which includes a positive definite matrix as parameters. In function approximation field, radial basis function (RBF) is widely used as basis. It performs well in unary and multivariate function approximation \cite{Light1992, WUZONGMIN1993Leef, Zhang1}. Compared to Gaussian density function, the number of parameters, instead of a positive definite matrix, reduces to only $1$. Thus when faced with high-dimensional problems, approximation using RBFs requires much fewer parameters than Gaussian density functions. Moreover, the method based on GMMs requires approximating both the two distributions with GMMs. However, using the basis functions in the same family to approximate two functions would not be reasonable in many situations. We may need to use two different families of basis functions to approximate the two functions. \cite{2024Chen} gives the optimal transport map and Wasserstein distance between distributions with radial density functions, which we call radially contoured distributions, and shows that the Wasserstein barycenter of such distributions is still radial. Thus we can establish a framework for finite mixtures of radially contoured distributions, which we call radial mixture models (RMMs).

In this paper, we focus on the space $W_2$. Instead of restricting the possible coupling measures to GMMs in \cite{DelonJulie2020Awdi}, we view RMMs as discrete distributions in the space of radially contoured distributions and define a distance for RMMs as the Wasserstein distance of their corresponding discrete distributions. We prove the parallel version for RMMs of important theorems for GMMs. Rather than requiring the components of identifiable mixtures of elliptically contoured distributions to be all of one type and satisfy marginal consistency in \cite{DelonJulie2020Awdi}, RMMs consist of components that can come from different families and do not need to satisfy marginal consistency. To infer the parameters of RMMs in application, we introduce the mini-batch stochastic technique to speed up the EM algorithm for RMMs. At last of this paper, we provide an application to color transfer task and compare our method with that proposed in \cite{DelonJulie2020Awdi} based on GMMs. The result shows that our method is more numerically stable and the transferred color distribution is more reasonable.

The paper is organized as follows. In \cref{sec2}, we recall the related results about Wasserstein distance and barycenter, especially for radially contoured distributions, and the Wasserstein-type distance $GW_2$ for GMMs proposed in \cite{DelonJulie2020Awdi}. In \cref{sec3}, we define $RW_2$ distance for identifiable finite mixtures of radially contoured distributions and prove a parallel version of the important properties of GMMs. We discuss in \cref{discussion} that we can remove the marginal consistency that is important in \cite{DelonJulie2020Awdi}. In \cref{sec4}, we conduct some toy experiments to show the difference between $RW_2$ and $W_2$. In \cref{sec5}, we show how to apply our method in practice. We compare our method with $GW_2$ in a color transfer task, and find that our method is more desirable in this experiment. Finally, we conclude our work in \cref{sec6} and provide some ways to which this work can be extended in the future.

\section{Review of related results}\label{sec2}

\subsection{Wasserstein distance and barycenter}

Let $d \geq 1$ be an integer. We write $\mathcal{P}(\R^d)$ the set of probability measures on $\R^d$. For $p \geq 1$, the Wasserstein space $\mathcal{P}_p(\R^d)$ is defined as the set of probability measures with finite moment of order $p$, that is,
\begin{equation*}
	\mathcal{P}_p(\R^d) = \left\{ \mu \in \mathcal{P}(\R^d): \int_{\R^d} \norm{x}^p \mathrm{d} \mu(x) < + \infty \right\}.
\end{equation*}
From now on, we will focus on the case $p=2$.

Given $\mu_0, \mu_1 \in \mathcal{P}_2(\R^d)$, a transport map $T$ from $\mu_0$ to $\mu_1$ is a measure preserving map, i.e. for every measurable set $A \subset \R_d$, $\mu_0(T^{-1}(A)) = \mu_1(A)$. The Monge formulation of optimal transport is
\begin{equation}\label{Monge}
	\inf_{T: T_\# \mu_0 = \mu_1} \int_{\R^d} \norm{x - T(x)}^2 \mathrm{d} \mu_0(x),
\end{equation}
where $T: T_\# \mu_0 = \mu_1$ means that $T$ takes over all measurable transport maps from $\mu_0$ to $\mu_1$. The optimal $T$ is called Monge map transporting $\mu_0$ to $\mu_1$. Define $\Pi(\mu_0,\mu_1) \subset \mathcal{P}_2(\R^d \times \R^d)$ as the subset of probability distributions $\gamma$ on $\R^d \times \R^d$ with marginals $\mu_0$ and $\mu_1$. Specifically, $(P_0)_\# \gamma = \mu_0$ and $(P_1)_\# \gamma = \mu_1$, where $P_0$ and $P_1$ are natural projections from $\R^d \times \R^d$ to $\R^d$. 
Kantorovich gives a relaxed form of (\ref{Monge}) by
\begin{equation}\label{W2}
	W_2(\mu_0, \mu_1)^2 = \inf_{\gamma \in \Pi(\mu_0, \mu_1)} \int_{\R^d\times\R^d} \norm{x-y}^2 \mathrm{d} \gamma(x,y).
\end{equation}
$W_2(\cdot,\cdot)$ actually defines a distance on $\mathcal{P}_2(\R^d)$, which is known as the Wasserstein distance.

If $T$ is a Monge map transporting $\mu_0$ to $\mu_1$, the path $(\mu_t)_{t\in[0,1]}$ given by
\begin{equation}\label{geodesic curve}
	\mu_t = (T_t)_\# \mu_0, \quad T_t(x) = (1-t)x + tT(x), \quad \forall \, t \in [0,1]
\end{equation}
defines a constant speed geodesic curve on $(\mathcal{P}_2(\R^d),W_2)$, i.e.
\begin{equation}\label{geodesic property}
	W_2(\mu_s, \mu_t) = (t-s) W_2(\mu_0,\mu_1), \quad 0 \leq s \leq t \leq 1.
\end{equation}
This path is called McCann interpolation \cite{McCannRobertJ.1997ACPf} between $\mu_0$ and $\mu_1$. It satisfies
\begin{equation}\label{mu_t}
	\mu_t \in \underset{\mu \in \mathcal{P}_2(\R^d)}{\arg\min} (1-t)W_2(\mu_0, \mu)^2 + tW_2(\mu_1,\mu)^2.
\end{equation}

The McCann interpolation can be generalized to more than two distributions, which is known as Wasserstein barycenter. Given an integer $N \geq 2$, $N$ probability measures $\mu_1, \mu_2, \cdots, \mu_N$ and weights $\lambda = (\lambda_1, \lambda_2, \cdots, \lambda_N) \in \Gamma_N$, where $\Gamma_N$ is the $N$-probability simplex, the associated Wasserstein barycenter problem is
\begin{equation}\label{Barycenter Problem}
	BW_2(\mu_1, \mu_2, \cdots, \mu_N)^2 = \inf_{\mu \in \mathcal{P}_2(\R^d)} \sum_{j=1}^N \lambda_j W_2(\mu_j,\mu)^2.
\end{equation}
A solution of the previous problem is called the Wasserstein barycenter of probability measures $\{\mu_j\}_{j=1}^N$ with weights $\{\lambda_j\}_{j=1}^N$. The existence and uniqueness have been deeply studied in \cite{AguehMartial2011Bitw}. 

\subsection{Results for radially contoured distributions}\label{subsec2.2}

The results in this section can be found in \cite{2024Chen}. A $d$-dimensional radially contoured distribution $\mu$ on $\R^d$ is a probability measure and
\begin{equation*}
	\mathrm{d}\mu(x) = \frac{1}{Z}\rho\left(\frac{\norm{x-m}}{c}\right) \dx,
\end{equation*}
where $\rho$ is a nonnegative measurable function or Dirac-$\delta$ function on $\R_+$ and $Z$ is the normalizer that ensures $\int_{\R^d} \mathrm{d} \mu(x) = 1$. In this paper, we mainly focus on the case that $\rho$ is continuous. For simplicity, we write $\rho(x;m,c)$ for $\rho(\norm{x-m}/c)$. We call $\rho$ the generator of $\mu$. We denote $\mu = R_d(m,c,\rho)$. If $c = 1$, we simply denote $\mu = R_d(m,\rho)$. Moreover, we will assume $\mu \in \mathcal{P}_2(\R^d)$. By straightforward calculation, we can get the normalizer, expectation and covariance of a radially contoured distributions as follows: let $X \sim \mu = R_d(m, c, \rho) \in \mathcal{P}_2(\R^d)$, then
\begin{equation*}
	Z = c^d |\mathbb{S}^{d-1}| \int_0^{+\infty} r^{d-1} \rho(r) \dr, \quad  E[X] = m, \quad \operatorname{Cov}[X] = \frac{c^2}{d} \frac{\int_0^{+\infty}r^{d+1}\rho(r)\dr}{\int_0^{+\infty}r^{d-1}\rho(r)\dr}I_d,
\end{equation*}
where $|\mathbb{S}^{d-1}| = 2 \pi^{d/2} / \Gamma(d/2)$ denotes the volume of $\mathbb{S}^{d-1}$.

Let $\mu_0 = R_d(m_0,\rho_0), \, \mu_1 = R_d(m_1,\rho_1)$ be two radially contoured distributions on $\R^d$. Considering two 1-dimensional radially contoured distributions $\tilde{\mu}_0 = R_1(0,r^{d-1}\rho_0(r))$ and $\tilde{\mu}_1 = R_1(0,r^{d-1}\rho_1(r))$, and the Monge map $C$ transporting $\tilde{\mu}_0$ to $\tilde{\mu}_1$, then the Monge map $T$ transporting $\mu_0$ to $\mu_1$ is given by
\begin{equation}\label{T}
	T(x) = \frac{C(\norm{x-m_0})}{\norm{x-m_0}}(x-m_0) + m_1.
\end{equation}
The McCann interpolation is given by $\mu_t = R_d(m_t,\rho_t)$, where $m_t = (1-t)m_0 + t m_1$ and $\rho_t$ satisfies
\begin{equation}\label{rho_t}
	\rho_t(C_t(r))C_t(r)^{d-1} \mathrm{d} C_t(r) = \rho_0(r) r^{d-1} \dr,
\end{equation}
where $C_t(r) = (1-t)r + tC(r)$. The Wasserstein distance is given by
\begin{equation*}
	W_2(\mu_0,\mu_1)^2 = \norm{m_0 - m_1}^2 + |\mathbb{S}^{d-1}| \int_0^{+\infty} (C(r) - r)^2 \rho_0(r) r^{d-1} \dr.
\end{equation*}
For $N \geq 2$, it has been proved in \cite{2024Chen} that the Wasserstein barycenter $\mu$ of radially contoured distributions $\{\mu_j = R_d(m_j,c_j,\rho_j)\}_{j=1}^N$ is still a radially contoured distribution centered at $m_* = \sum_{j=1}^N \lambda_j m_j$.

When the generators of radially contoured distributions are the same, the conclusions above have a simple form. Let $\mu_0 = R_d(m_0, c_0, \rho)$, $\mu_1 = R_d(m_1, c_1, \rho)$, then the Monge map from $\mu_0$ to $\mu_1$ is
\begin{equation*}
	T(x) = \frac{c_1}{c_0}(x-m_0) + m_1.
\end{equation*}
The displacement interpolation is $\mu_t = R_d(m_t,c_t,\rho)$, where $c_t = (1-t) m_0 + t m_1$, and
\begin{equation*}
	W_2(\mu_0,\mu_1)^2 = \norm{m_0 - m_1}^2 + (c_0 - c_1)^2 \frac{\int_0^{+\infty} r^{d+1} \rho(r)\dr}{\int_0^{+\infty} r^{d-1} \rho(r) \dr}.
\end{equation*}
For $N \geq 2$, let $\mu_j = R_d(m_j,c_j,\rho_j)$, $j = 1, 2, \cdots, N$, and $\lambda = (\lambda_1, \lambda_2, \cdots, \lambda_N)$ be the weights. There exists $\rho$ such that barycenter of $\{\mu_j\}_{j=1}^N$ with weights $\{\lambda_j\}_{j=1}^N$ is $\mu = R_d(m_*,c_*,\rho)$, where
\begin{equation*}
	m_* = \sum_{j=1}^N \lambda_j m_j, \quad c_* = \sum_{j=1}^N \lambda_j c_j.
\end{equation*}
When $\rho_1 = \rho_2 = \cdots = \rho_N$, the barycenter is simply $R_d(m_*, c_*, \rho_1)$.

\subsection{Wasserstein-type distance for GMM}

Let $g_{m, \Sigma}$ be the density function of Gaussian distribution with mean $m$ and covariance $\Sigma$. A Gaussian mixture model (GMM) of $K$ components is a probability distribution with density $\rho(x) = \sum_{i=1}^K \pi_i g_{m_i, \Sigma_i}(x)$, where $\pi = (\pi_1, \pi_2, \cdots, \pi_K)$ is a weight vector. Let $\mu_0 = \sum_{k=1}^{K_0} \pi_0^k \nu_0^k$ and $\mu_1 = \sum_{l=1}^{K_1} \pi_1^l \nu_1^l$ be two $d$-dimensional GMMs. In \cite{DelonJulie2020Awdi}, the authors have proposed a Wasserstein-type distance $MW_2$, which is formulated as
\begin{equation}\label{MW2}
	MW_2(\mu_0, \mu_1)^2 = \inf_{\gamma \in \Pi(\mu_0, \mu_1) \cap GMM_{2d}(\infty)} \int_{\R^d \times \R^d} \norm{x - y}^2 \mathrm{d} \gamma(x,y),
\end{equation}
where $\mu_0, \mu_1$ are $d$-dimensional Gaussian mixture models and $GMM_d(\infty)$ denotes the set of finite mixtures of $d$-dimensional Gaussian distributions. They show that (\ref{MW2}) has an equivalent discrete formulation: we can view a GMM as a discrete measure in the space of Gaussian distributions and compute the discrete optimal transport, i.e.
\begin{equation}\label{MW2 discrete}
	MW_2(\mu_0, \mu_1)^2 = \inf_{w \in \Pi(\pi_0, \pi_1)} \sum_{k,l = 1}^{K_0, K_1} w_{kl} W_2(\nu_0^k, \nu_1^l)^2.
\end{equation}
$MW_2$ actually defines a metric on $GMM_d(\infty)$, and $GMM_d(\infty)$ equipped with $MW_2$ is a geodesic space. The natural $MW_2$ barycenter problem is
\begin{equation}\label{MW barycenter}
	\inf_{\mu \in GMM_d(\infty)} \sum_{j=1}^N \lambda_j MW_2(\mu_j,\mu)^2,
\end{equation}
and the solution of (\ref{MW barycenter}) is called the $MW_2$ barycenter for $N$ GMMs $\{\mu_j\}_{j=1}^N$, where $\mu_j = \sum_{j=1}^{K_j} \pi_j^{k_j} \nu_j^k$, with weights $\{\lambda_1, \lambda_2, \cdots, \lambda_N\}$. They give a discrete formulation of (\ref{MW barycenter}) which can be computed via linear programming:
\begin{equation*}
	\inf_{\mu \in GMM_d(\infty)} \sum_{j=1}^N \lambda_j MW_2(\mu_j,\mu)^2 = \min_{w\in\Pi(\pi_1,\pi_2,\cdots,\pi_N)} \sum_{k_1, \cdots, k_N = 1}^{K_1, \cdots, K_N} w_{k_1 \cdots k_N} BW_2(\nu_1^{k_1}, \cdots, \nu_N^{k_N})^2.
\end{equation*}
They also consider about how to define a map via $MW_2$. Let $\gamma^*$ and $w^*$ be the optimal solutions of \cref{MW2} and \cref{MW2 discrete}, and $T_{k,l}$ be the Monge map transporting $\nu_0^{k}$ to $\nu_1^l$. Define a discrete probability distribution $\{p_{k,l}(x)\}$ as
\begin{equation}\label{pkl GMM}
	p_{k,l}(x) = \frac{w^*_{k,l} g_{m_0^k, \Sigma_0^k} (x)}{\sum_{j} \pi_0^j g_{m_0^k, \Sigma_0^k} (x)}.
\end{equation}
They construct two kinds of maps. The first one is $T_{\rm mean}$:
\begin{equation*}
	T_{\rm mean}(x) = \sum_{k,l = 1}^{K_0, K_1} p_{k,l}(x) T_{k,l}(x),
\end{equation*}
and the second one is $T_{\rm rand}$:
\begin{equation*}
	T_{\rm rand}(x) = T_{k,l}(x) \quad \text{with probability} \quad p_{k,l}(x).
\end{equation*}
Compared to $T_{\rm mean}$, $T_{\rm rand}$ is more mathematically reasonable. They show that for every measurable set $A$ of $\R^d$,
\begin{equation*}
	\int_{\R^d} P(T_{\rm rand}(x) \in A) \mathrm{d} \mu_0(x) = \mu_1(A).
\end{equation*}
It means that $({T_{\rm rand}})_\# \mu_0$ is equal to $\mu_1$ in expectation. But in practice, $T_{\rm rand}$ gives a noisier result than $T_{\rm mean}$. In \cite{DelonJulie2020Awdi}, $T_{\rm mean}$ is used in most of the experiments.

\section{The space of radial mixture models equipped with $RW_2$ distance}\label{sec3}

\subsection{The definition of radial mixture models and $RW_2$ distance}

\begin{definition}
	Let $K \geq 1$ be an integer. A finite radial mixture $\mu$ of size $K$ on $\R^d$ is an identifiable probability measure that can be written as
	\begin{equation*}
		\mathrm{d} \mu(x; m, c,\rho) = \sum_{k=1}^K \pi^k \mathrm{d} \nu^k(x; m^k, c^k,\rho^k),
	\end{equation*}
	where $\nu^k = R_d(m^k,c^k,\rho)$, $\pi \in \Gamma_K$, and $\rho$ is a Dirac $\delta$ function or nonnegative continuous function which ensures $\mu \in \mathcal{P}_2(\R^d)$. Henceforth, in certain situations, parameters $m, c$ and $\rho$ of $\mu$ or $\nu$ can be omitted. We also call $\rho$ the generator of radial mixture $\mu$. We denote the set of all radial mixtures generated from $\rho$ on $\R^d$ as $RMM_d(\rho)$ and
	\begin{equation*}
		RMM_d = \cup_{\rho} RMM_d(\rho).
	\end{equation*}
	Here $\rho$ takes over all functions that promise $RMM_d(\rho)$ is identifiable, i.e. given two mixtures $\mu_0$ and $\mu_1$ in $RMM_d(\rho)$, where
	\begin{equation*}
		\mathrm{d}\mu_j(x) = \sum_{k=1}^{K_j} \pi_j^k \mathrm{d} \nu_j^k(x), \quad \nu_j^k = R_d(m_j^k,c_j^k,\rho), \quad j = 0,1,
	\end{equation*}
	and $\nu_j^{k} \neq \nu_j^l$ for all $1 \leq k < l \leq K_j$, then $\mu_0 = \mu_1$ if and only if $K_0 = K_1$ and there exists a $K_0$-permutation $\sigma$ such that $\pi_0^{\sigma(k)} = \pi_1^k$ and $\nu_0^{\sigma(k)} = \nu_1^k$ for all $1 \leq k \leq K_0$.
\end{definition}

There exist many choices of such $\rho$. We can view radially contoured distributions as a degeneration of elliptically contoured distributions. The identifiability of finite mixtures of elliptically contoured distributions has been discussed in \cite{HOLZMANNHAJO2006IoFM}. For instance, finite mixtures of multi-variate t-distributions are identifiable. The following proposition allows us to approximate distributions in $\mathcal{P}_2(\R^d)$ with radial mixtures.
\begin{proposition}\label{thm:dense}
	$RMM_d(\rho)$ is dense in $\mathcal{P}_2(\R^d)$ with the distance $W_2$.
\end{proposition}

The proof is similar to that of Proposition $1$ in \cite{DelonJulie2020Awdi}. By viewing radial mixtures as discrete distributions on the space of radially contoured distributions, we deduce a distance for radial mixtures. This distance is analogous to the Wasserstein distance for discrete distributions on Euclidean space.
\begin{theorem}\label{RW_2}
	Let $\mu_0, \mu_1 \in RMM_d$, where
	\begin{equation}\label{mui}
		\mathrm{d}\mu_j(x) = \sum_{k=1}^{K_j} \pi_j^k \mathrm{d} \nu_j^k(x), \quad \nu_j^k = R_d(m_j^k,c_j^k,\rho_j), \quad j = 0,1,
	\end{equation}
	and $w^*$ is a solution of
	\begin{equation}\label{w*}
		\min_{w\in\Pi(\pi_0,\pi_1)} \sum_{k,l = 1}^{K_0,K_1} w_{kl} W_2(\nu_0^k,\nu_1^l)^2,
	\end{equation}
	then
	\begin{equation*}
		RW_2(\mu_0,\mu_1) = \sqrt{\sum_{k,l=1}^{K_0,K_1} w^*_{kl} W_2(\nu_0^k,\nu_1^l)^2}
	\end{equation*}
	is a distance on $RMM_d$. We call $w \in \Pi(\pi_0,\pi_1)$ the $RW_2$ plan between mixtures $\mu_0$ and $\mu_1$ and $w^*$ the $RW_2$ optimal plan.
\end{theorem}

\begin{proof}
	Apparently, $RW_2(\mu_0,\mu_1) = RW_2(\mu_1, \mu_0)$. For the positive definite property, we will first prove $\mu_0 = \mu_1$ if $RW_2(\mu_0,\mu_1) = 0$. Let $T_{kl}$ be the Monge map from $\nu_0^k$ to $\nu_1^l$. We define
	\begin{equation}\label{gamma}
		\mathrm{d} \gamma(x,y) = \sum_{k,l} w^*_{kl} \delta(y - T_{kl}(x)) \mathrm{d} \nu_0^k(x) \dy ,
	\end{equation}
	then for all $h \in C(\R^d)$, we have
	\begin{equation*}
		\begin{aligned}
			& \int_{\R^d \times \R^d} h(x) \mathrm{d} (P_0)_{\#}\gamma(x,y)   = \sum_{k,l} \int_{\R^d\times\R^d} h(x) w^*_{kl} \delta(y - T_{kl}(x)) \mathrm{d} \nu_0^k(x) \dy \\ 
			= & \sum_{k,l} \int_{\R^d} h(x) w^*_{kl} \mathrm{d} \nu_0^k(x) = \sum_k \int_{\R^d} h(x) \pi_0^k \mathrm{d} \nu_0^k(x) = \int_{\R^d} h(x) \mathrm{d} \mu_0(x), \\
			& \int_{\R^d \times \R^d} h(y) \mathrm{d} (P_1)_{\#}\gamma(x,y)   = \sum_{k,l} \int_{\R^d\times\R^d} h(y) w^*_{kl} \delta(y - T_{kl}(x)) \mathrm{d} \nu_0^k(x) \dy \\ 
			= & \sum_{k,l} \int_{\R^d} h(y) w^*_{kl} \mathrm{d} \nu_1^l(y) = \sum_l \int_{\R^d} h(y) \pi_1^l \mathrm{d} \nu_1^l(y) = \int_{\R^d} h(y) \mathrm{d} \mu_1(y).
		\end{aligned}
	\end{equation*}
	Therefore, $\gamma \in \Pi(\mu_0,\mu_1)$. Moreover,
	\begin{equation*}
		\begin{aligned}
			& \int_{\R^d\times\R^d} \norm{x-y}^2 \mathrm{d} \gamma(x,y) = \sum_{k,l} w^*_{kl} \int_{\R^d\times\R^d} \norm{x-y}^2 \delta(y - T_{kl}(x)) \mathrm{d} \nu_0^k(x) \dy \\ 
			= & \sum_{k,l} w^*_{kl} \int_{\R^d} \norm{x - T_{kl}(x)}^2 \mathrm{d} \nu_0^k(x) = \sum_{k,l} w^*_{kl} W_2(\nu_0^k,\nu_1^l)^2 = RW_2(\mu_0,\mu_1)^2.
		\end{aligned}  
	\end{equation*}
	It follows that
	\begin{equation*}
		W_2(\mu_0,\mu_1)^2 \leq \int_{\R^d\times\R^d} \norm{x-y}^2 \mathrm{d} \gamma(x,y) = RW_2(\mu_0,\mu_1)^2.
	\end{equation*}
	Thus if $RW_2(\mu_0,\mu_1) = 0$, we have $W_2(\mu_0,\mu_1) = 0$, and thus $\mu_0 = \mu_1$. On the other hand, if $\mu_0 = \mu_1$, by the identifiablility of $RMM_d$, after resorting the index, we can assume $K_0 = K_1$, $\pi_0 = \pi_1$ and $\nu_0^k = \nu_1^k$ for all $k$. Then we can choose $w$ as the identity matrix, which means $RW_2(\mu_0, \mu_1) \leq 0$.
	
	Lastly, we only need to prove the triangle inequality, i.e.
	\begin{equation*}
		RW_2(\mu_0, \mu_1) + RW_2(\mu_1,\mu_2) \geq RW_2(\mu_0, \mu_2)
	\end{equation*}
	for any $\mu_0, \mu_1, \mu_2 \in RMM_d$. Denote the probability vector associated with  $\mu_0,\mu_1,\mu_2$ by $\pi_0,\pi_1,\pi_2$, the components of $\mu_j$ by $\nu_j^k$ and the sizes by $K_0, K_1, K_2$. Let $w^{01}(w^{12})$ be the solution to (\ref{w*}) with marginals $\pi_0, \pi_1 (\pi_1, \pi_2)$. Define $w^{02}$ by
	\begin{equation*}
		w^{02}_{kl} = \sum_{j=1}^{K_1} \frac{w_{kj}^{01} w_{jl}^{12}}{\pi_1^j}.
	\end{equation*}
	If $\pi_1^j = 0$, the corresponding term is set to $0$. Then $w^{02} \in \Pi(\pi_0,\pi_2)$, which follows from
	\begin{equation*}
		\sum_{k} w_{kl}^{02} = \sum_{j} w_{jl}^{12} = \pi_2^l, \quad \forall \, l = 1, 2, \cdots, K_2,
	\end{equation*}
	\begin{equation*}
		\sum_{l} w_{kl}^{02} = \sum_{j} w_{kj}^{01} = \pi_0^k. \quad \forall \, k = 1, 2, \cdots, K_0.
	\end{equation*}
	Therefore,
	\begin{align*}
		RW_2(\mu_0,\mu_2) & \leq \sqrt{\sum_{k,l} w^{02}_{kl} W_2(\nu_0^k,\nu_2^l)^2} = \sqrt{\sum_{k,l,j} \frac{w_{kj}^{01} w_{jl}^{12}}{\pi_1^j} W_2(\nu_0^k,\nu_2^2)^2} \\
		& \leq \sqrt{\sum_{k,l,j} \frac{w_{kj}^{01} w_{jl}^{12}}{\pi_1^j} \left[ W_2(\nu_0^k, \nu_1^j) + W_2(\nu_1^j, \nu_2^l) \right]^2} \\
		& \leq \sqrt{\sum_{k,l,j} \frac{w_{kj}^{01} w_{jl}^{12}}{\pi_1^j} W_2(\nu_0^k,\nu_1^j)^2 } + \sqrt{\sum_{k,l,j} \frac{w_{kj}^{01} w_{jl}^{12}}{\pi_1^j} W_2(\nu_1^j,\nu_2^l)^2} \\
		& = \sqrt{\sum_{k,j}w_{kj}^{01} W_2(\nu_0^k,\nu_1^j)^2} + \sqrt{\sum_{l,j} w_{jl}^{12} W_2(\nu_1^j,\nu_2^l)^2 } \\
		& = RW_2(\mu_0, \mu_1) + RW_2(\mu_1, \mu_2).
	\end{align*}
	In the above, the second inequality is due to the triangle inequality of $W_2$, and the third inequality comes from the Minkowski inequality.
\end{proof}

Note that this distance is quite easier to compute than Wasserstein distance, since we do not need to discretize the radial mixtures in space. We only need to compute the Monge transport map for two $1$-dimensional distributions as discussed in \cref{subsec2.2} and a linear programming with $K_0 \cdot K_1$ variables.


    \subsection{The space $(RMM_d, RW_2)$ is a geodesic space}

\begin{theorem}\label{RW2 geodesic}
	Let $\mu_0, \mu_1 \in RMM_d$ be written as in (\ref{mui}). The displacement interpolation on $(RMM_d, RW_2)$ connecting $\mu_0$ and $\mu_1$ is given by
	\begin{equation}\label{RW2 interpolation}
		\mu_t = \sum_{k,l} w^*_{kl} \nu_t^{kl},
	\end{equation}
	where $\nu_t^{kl}$ is the displacement interpolation between $\nu_0^k$ and $\nu_1^l$, and
	\begin{equation}\label{thm2}
		RW_2(\mu_s,\mu_t) = (t-s) RW_2(\mu_0, \mu_1), \quad 0 \leq s \leq t \leq 1.
	\end{equation}
\end{theorem}

\begin{proof}
	For any $0 \leq s \leq t \leq 1$, we have
	\begin{equation*}
		RW_2(\mu_s,\mu_t) \leq \sqrt{\sum_{k,l} w^*_{kl} W_2(\nu_s^{kl},\nu_t^{kl})^2} = (t-s) \sqrt{\sum_{k,l} w^*_{kl}W_2(\nu_0^k,\nu_1^l)^2} = (t-s) RW_2(\mu_0,\mu_1).
	\end{equation*}
	The inequality is given by choosing the transport plan as transporting all mass of $\nu_s^{kl}$ to $\nu_t^{kl}$, and the equality is given by (\ref{geodesic property}). It follows that
	\begin{equation*}
		\begin{aligned}
			& RW_2(\mu_0,\mu_s) + RW_2(\mu_s, \mu_t) + RW_2(\mu_t, \mu_1) \\ \leq & sRW_2(\mu_0,\mu_1) + (t-s) RW_2(\mu_0, \mu_1) + (1-t) RW_2(\mu_0, \mu_1) = RW_2(\mu_0,\mu_1).
		\end{aligned}
	\end{equation*}
	On the other hand, according to (\ref{RW_2}), we have
	\begin{equation*}
		RW_2(\mu_0,\mu_s) + RW_2(\mu_s,\mu_t) + RW_2(\mu_t,\mu_1) \geq RW_2(\mu_0,\mu_1).
	\end{equation*}
	Therefore, we obtain (\ref{thm2}).
\end{proof}


Compared to the geodesic path defined by $W_2$, the distributions on the $RW_2$ geodesic all keep the radial mixture structure. The examples can be seen in \cref{sec4}.

\subsection{The transport map defined in $(RMM_d, RW_2)$}\label{subsec33}

In many cases, we not only need the transport plan, but also need a transport map. Let $\mu_0$ and $\mu_1$ be two radial mixtures. We have proved that $\gamma$ defined in (\ref{gamma}) is a transport plan of the two RMM marginals. Thus we have two natural ways to define the transport map. The first way is setting
\begin{equation}
	T_{\rm mean}(x) = \mathbb{E}_\gamma (Y|X=x),
\end{equation}
where $(X,Y)$ is distributed according to $\gamma$ in (\ref{gamma}). The distribution of $Y|X=x$ is given by the discrete distribution
\begin{equation}\label{prob}
	\sum_{k,l} p_{k,l}(x) \delta_{T_{kl}(x)} \quad \text{with} \quad p_{k,l}(x) = \frac{\sum_{k,l} w_{k,l}^* \rho_0(x;m_0^k,c_0^k)}{\sum_j \pi_0^j \rho_0(x;m_0^k,c_0^k)}.
\end{equation}
We can get that
\begin{equation}\label{Tmean}
	T_{\rm mean}(x) = \frac{\sum_{k,l} w_{k,l}^* \rho_0(x;m_0^k,c_0^k) T_{kl}(x)}{\sum_j \pi_0^j \rho_0(x;m_0^k,c_0^k)}.
\end{equation}
This seems to be the most natural way to define a transport map. However, $({T_{\rm mean}})_\# \mu_0$ may be far from $\mu_1$. For instance, given three radially contoured distributions $\nu_0 = R_1(0,1,\rho), \nu_1 = R_1(a,1,\rho)$ ,$ \nu_2 = R_1(-a,1,\rho)$, we define $\mu_0 = \nu_0$ and $\mu_1 = 0.5 \nu_1 + 0.5 \nu_2$. In this case, $T_{\rm mean}$ is the identity map, and thus $({T_{\rm mean}})_\# \mu_0 = \mu_0$, which will be very far from $\mu_1$ when $a$ is large.

Based on (\ref{prob}), another way to define a transport map is setting
\begin{equation}\label{Trand}
	T_{\rm rand}(x) = T_{kl}(x) \quad \text{with probability} \, \, \, p_{k,l}(x) = \frac{\sum_{k,l} w_{k,l}^* \rho_0(x;m_0^k,c_0^k)}{\sum_j \pi_0^j \rho_0(x;m_0^k,c_0^k)}.
\end{equation}
As it has been discussed in \cite{DelonJulie2020Awdi}, $(T_{\rm rand})_\# \mu_0$ would be equal in expectation to $\mu_1$.

\subsection{The comparison between $RW_2$ and $W_2$}

\begin{proposition}[The comparison between $RW_2$ and $W_2$]
	Let $\mu_0, \mu_1 \in RMM_d$ be written as in (\ref{mui}), then
	\begin{equation}
			W_2(\mu_0,\mu_1) \leq RW_2(\mu_0,\mu_1) \leq W_2(\mu_0,\mu_1) + \sum_{j=0,1}\left[ 2\sum_{k=1}^{K_j} \pi_j^k (c_j^k)^2 \frac{\int_0^{+\infty} r^{d+1}\rho_j(r)\dr}{\int_0^{+\infty} r^{d-1}\rho_j(r)\dr} \right]^\frac{1}{2}.
		\end{equation}
\end{proposition}
\begin{proof}
	The first inequality has already been proved in (\ref{RW_2}). The proof of the second inequality is similar to the proof of Proposition $6$ in \cite{DelonJulie2020Awdi}.
\end{proof}

The first inequality shows that $RW_2$ distance is lower bounded by $W_2$ distance. The second inequality shows that when given two generators, $RW_2$ distance is upper bounded by $W_2$ plus a term that only depends on the parameters $\{c_j^k\}$. Therefore, when $\max_{j,k} \{c_j^k\} \leq \varepsilon$, we have
\begin{equation*}
	RW_2(\mu_0,\mu_1) \leq W_2(\mu_0,\mu_1) + C \varepsilon,
\end{equation*}
where $C$ is a constant depending on the dimension $d$ and generators $\rho_0, \rho_1$.

\subsection{The barycenter in $(RMM_d, RW_2)$}

\begin{definition}[$RW_2$ barycenter]
	Let $\mu_1, \mu_2, \cdots, \mu_N \in RMM_d$, where
	\begin{equation*}
		\mu_j = \sum_{k=1}^{K_j} \pi_j^k \nu_j^k, \quad \nu_j^k = R_d(m_j^k, c_j^k, \rho_j), \quad j = 1, 2, \cdots, N,
	\end{equation*}
	and $\lambda = (\lambda_1, \lambda_2, \cdots, \lambda_N) \in \Gamma_K$ is a weight vector. The $RW_2$ barycenter $\mu$ of the $\mu_j$ with weights $\{\lambda_j\}$ is a solution of
	\begin{equation}\label{RW2 barycenter}
		\inf_{\mu \in RMM_d} \sum_{j=1}^N \lambda_j RW_2(\mu_j, \mu)^2.
	\end{equation}
\end{definition}

The definition of $RW_2$ barycenter involves a variational problem, wich is difficult to handle. To get a more specific description of the $RW_2$ barycenter, we consider a multimarginal problem as follows:
\begin{equation}\label{MRW2}
	MRW_2(\mu_1, \mu_2, \cdots, \mu_N)^2 = \min_{w\in\Pi(\pi_1,\pi_2,\cdots,\pi_N)} \sum_{k_1, \cdots, k_N = 1}^{K_1, \cdots, K_N} w_{k_1 \cdots k_N} BW_2(\nu_1^{k_1}, \cdots, \nu_N^{k_N})^2,
\end{equation}
where $\Pi(\pi_1,\pi_2,\cdots,\pi_N)$ is the set of tensors $w$ having $\pi_1, \pi_2, \cdots, \pi_N$ as marginals, i.e.
\begin{equation*}
	\sum_{\substack{k_1,\cdots,k_{j-1},k_{j+1},\cdots,k_N \\ k_j = k}} w_{k_1 \cdots k_N} = \pi_j^k, \quad \forall \, j \in \{1, \cdots, N\}, k \in \{1, \cdots, K_j\}.
\end{equation*}

The two problems are linked by the following theorem.

\begin{theorem}\label{thm:RW2 barycenter}
	Let $\mu_1, \mu_2, \cdots, \mu_N \in RMM_d$, then
	\begin{equation*}
		\inf_{\mu \in RMM_d} \sum_{j=1}^N \lambda_j RW_2(\mu_j, \mu)^2 = MRW_2(\mu_1, \mu_2, \cdots, \mu_N)^2.
	\end{equation*}
	Let $w^* \in \Pi(\pi_1,\cdots,\pi_N)$ be the solution of (\ref{MRW2}), and $\nu^{k_1 \cdots k_N}$ be the barycenter of $\nu_1^{k_1}, \cdots, \nu_N^{k_N}$ with weights $\{\lambda_j\}_{j=1}^N$. The $RW_2$ barycenter $\mu$ of $\mu_j$ with weights $\{\lambda_j\}_{j=1}^N$ is given by
	\begin{equation*}
		\mu = \sum_{k_1 \cdots k_N} w_{k_1 \cdots k_N}^* \nu^{k_1 \cdots k_N}.
	\end{equation*}
	Moreover, the number of components of this barycenter is no more than $K_1 + K_2 + \cdots + K_N - N + 1$.
\end{theorem}

\begin{proof}
	For any $\mu \in RMM_d$, we write $\mu = \sum_{l=1}^L \pi^l \nu^l$. Let $w^j$ be the $RW_2$ optimal plan between the mixtures $\mu_j$ and $\mu$. We define a $K_1 \times \cdots \times K_N \times L$ tensor $\alpha$ and a $K_1 \times \cdots \times K_N$ tensor $\bar{\alpha}$ by
	\begin{equation*}
		\alpha_{k_1 \cdots k_N l} = \frac{\Pi_{j=1}^{N} w_{k_j l}^{j}}{(\pi^l)^N}, \quad \text{and} \quad \bar{\alpha}_{k_1\cdots k_N} = \sum_{l=1}^L \alpha_{k_1 \cdots k_N l}.
	\end{equation*}
	Clearly $\alpha \in \Pi(\pi_1,\cdots,\pi_N,\pi)$ and $\bar{\alpha} \in \Pi(\pi_1,\cdots,\pi_N)$. Moreover,
	\begin{equation*}
		\begin{aligned}
			& \sum_{j=1}^N \lambda_j RW_2(\mu_j,\mu)^2 = \sum_{j=1}^N \lambda_j \sum_{k_j=1}^{K_j} \sum_{l=1}^L w_{k_j l}^j W_2(\nu_j^{k_j},\nu^l)^2 \\
			= & \sum_{j=1}^{N} \lambda_j \sum_{k_1 \cdots k_N,l} \alpha_{k_1 \cdots k_N l} W_2(\nu_j^{k_j},\nu^l)^2 = \sum_{k_1 \cdots k_N,l} \alpha_{k_1 \cdots k_N l} \sum_{j=1}^N \lambda_j W_2(\nu_j^{k_j},\nu^l)^2 \\
			\geq & \sum_{k_1 \cdots k_N,l} \alpha_{k_1 \cdots k_N l} BW_2(\nu_1^{k_1},\cdots,\nu_N^{k_N})^2 = \sum_{k_1 \cdots k_N} \bar{\alpha}_{k_1 \cdots k_N} BW_2(\nu_1^{k_1},\cdots,\nu_N^{k_N})^2 \\
			\geq & MRW_2(\mu_1,\cdots,\mu_N)^2.
		\end{aligned}
	\end{equation*}
	This inequality holds for all $\mu \in RMM_d$. Thus
	\begin{equation*}
		\inf_{\mu \in RMM_d} \sum_{j=1}^N \lambda_j RW_2(\mu_j, \mu)^2 \geq MRW_2(\mu_1, \mu_2, \cdots, \mu_N)^2.
	\end{equation*}
	Conversely, let $w^* \in \Pi(\pi_1,\cdots,\pi_N)$ be the solution of (\ref{MRW2}) and $\nu^{k_1 \cdots k_N}$ be the barycenter of $\nu_1^{k_1}, \cdots, \nu_N^{k_N}$ with weights $\{\lambda_j\}_{j=1}^N$. We define $\mu = \sum_{k_1 \cdots k_N} w_{k_1 \cdots k_N}^* \nu^{k_1 \cdots k_N}$. Based on the results in \cref{subsec2.2}, $\nu^{k_1 \cdots k_N}$ is still a radially contoured distribution, and thus $\mu \in RMM_d$. We define $RW_2$ plans $\{w^j\}_{j=1}^N$ between $\{\mu_j\}_{j=1}^N$ and $\mu$ by
	
	\begin{equation*}
		w^j_{i_j,k_1\cdots k_N} =\begin{cases}
			w^*_{k_1 \cdots k_N},& \quad i_j = k_j, \\
			0,& \quad i_j \neq k_j,
		\end{cases}
	\end{equation*}
	where $1 \leq i_j \leq K_j$. It is clearly $w^j \in \Pi(\pi_j, w)$. Then we have
	
	\begin{equation*}
		\begin{aligned}
			& MRW_2(\mu_1, \mu_2, \cdots, \mu_N)^2 =  \sum_{k_1 \cdots k_N} w_{k_1 \cdots k_N}^* BW_2(\nu_1^{k_1},\cdots,\nu_N^{k_N})^2\\
			= & \sum_{k_1 \cdots k_N} w_{k_1 \cdots k_N}^* \sum_{j=1}^N \lambda_j W_2(\nu_j^{k_j},\nu^{k_1 \cdots k_N})^2 = \sum_{j=1}^N \lambda_j \sum_{k_1 \cdots k_N} w_{k_1 \cdots k_N}^* W_2(\nu_j^{k_j},\nu^{k_1 \cdots k_N})^2 \\
			= & \sum_{j=1}^N \lambda_j \sum_{i_j}\sum_{k_1 \cdots k_N} w^j_{i_j,k_1 \cdots k_N} W_2(\nu_j^{i_j},\nu^{k_1 \cdots k_N})^2 \geq \sum_{j=1}^N \lambda_j RW_2(\mu_j,\mu)^2 \\ 
			\geq &0 \inf_{\mu \in RMM_d} \sum_{j=1}^N \lambda_j RW_2(\mu_j,\mu)^2.
		\end{aligned}
	\end{equation*}
	Therefore, 
	\begin{equation*}
		\inf_{\mu \in RMM_d} \sum_{j=1}^N \lambda_j RW_2(\mu_j,\mu)^2 = MRW_2(\mu_1,\mu_2,\cdots,\mu_N)^2.
	\end{equation*}
	Moreover, since $MRW_2(\mu_1,\mu_2,\cdots,\mu_N)^2 \geq \sum_{j=1}^N \lambda_j RW_2(\mu_j, \mu)^2$, we get that $\mu$ is a $RW_2$ barycenter of $\mu_j$ with weights $\{\lambda_j\}_{j=1}^N$. Since (\ref{MRW2}) is a linear program with $K_1 + \cdots + K_N - N + 1$ linear constraints, the solution $w^*$ can not have more than  $K_1 + \cdots + K_N - N + 1$ nonzero elements. Therefore, the barycenter has more than $K_1 + \cdots + K_N - N + 1$ components.
\end{proof}

\cref{thm:RW2 barycenter} shows that the $RW_2$ barycenter is a combination of the Wasserstein barycenters of the components, and the functional optimal problem \cref{RW2 barycenter} can be reduced to a linear programming problem. Therefore, to compute an $RW_2$ barycenter, we only need to compute the Wasserstein barycenters of radial components of the marginals and solve a linear programming problem. When the marginals come from the same family $RMM_d(\rho)$, their $RW_2$ barycenter can be easily computed since the Wasserstein barycenter of radially contoured distributions $\{\mu_j = R_d(m_j, c_j, \rho)\}_{j=1}^N$ with weights $\{\lambda_j\}_{j=1}^N$ is simply $R_d(m_*,c_*,\rho)$, where $m_* = \sum_{j=1}^N \lambda_j m_j$ and $c_* = \sum_{j=1}^N \lambda_j c_j$. In this case, only a linear programming problem need to be solved, and thus the $RW_2$ barycenter is quite easier to compute than Wasserstein barycenter. In addition, the Wasserstein barycenter of radial mixtures is usually not a radial mixtures. However, our $RW_2$ will keep this structure.

\subsection{Some remarks on our method}\label{discussion}

In the framework given in \cite{DelonJulie2020Awdi}, the mixtures are identifiable and the components satisfy marginal consistency, where the marginal consistency requirement is important in the $GW_2$ definition. Their definition of $GW_2$ between two $d$-dimensional mixtures from the same family restricts the transport plan to the space of $2d$-dimensional mixtures. However, our $RW_2$ is not defined in that way. In our framework, only the identifiability is required. To show the importance of identifiability of $RMM_d(\rho)$, we consider the following instance. Let $\rho_{U}(x) = \mathbbm{1}_{[0,1)}(x)$ be the indicator function on $[0,1)$ and let $\mu_0$ be the uniform distribution on $(-1,1)$, which can be written as $R_1(0,1,\rho)$. Define $\nu^0 = R_1(0,2,\rho_U)$, $\nu^1 = R_1(1,2,\rho_U)$, and
\begin{equation*}
	\mathrm{d}\mu_1(x) = 0.5 \mathrm{d} \nu^0(x) + 0.5 \mathrm{d} \nu^1(x).
\end{equation*}
We then have $\mu_1 \in RMM_1(\rho_U)$ and $\mu_0 = \mu_1$, which means that the mixtures in $RMM_1(\rho_U)$ are not identifiable. However, by the definition of $RW_2$ in \cref{w*}, $w^*$ can only be chosen as $(0.5, 0.5)$, which follows that
\begin{equation*}
	RW_2(\mu_0,\mu_1) = \sqrt{0.5 W_2(\mu_0, \nu^0)^2 + 0.5 W_2(\mu_0, \nu^1)^2} > 0,
\end{equation*}
and thus $RW_2$ is not a distance on $RMM_1(\rho_U)$.

Another question is whether we can generalize $RW_2$ from the radially contoured case to the general elliptically contoured case. An elliptically contoured distribution on $\R^d$ has a density of the form
\begin{equation*}
	f_{m,\Sigma}(x) = \frac{1}{Z} \rho\left(\left(x-m\right)^\top \Sigma^{-1} \left(x-m\right)\right), \quad \forall x \in \R^d,
\end{equation*}
where $Z$ is the normalizer, $m \in \R^d$, $\Sigma$ is a positive definite matrix of size $d \times d$, and $\rho$ is a positive generator defined on $\R_+$. We denote the set of all identifiable finite mixtures of elliptically contoured distributions with generator $\rho$ as $EMM_d(\rho)$ and the union of such $EMM_d(\rho)$ as $EMM_d$. For example, when choosing $\rho(x) = \exp(-x)$, $EMM_d(\rho) = GMM_d$. The generalization to $EMM_d(\rho)$ is almost trivial thanks to the results in \cite{GelbrichMatthias1990OaFf, GomezEusebio2003Asoc, AguehMartial2011Bitw}. \cite{GelbrichMatthias1990OaFf} and \cite{GomezEusebio2003Asoc} give the distance and the Monge map of two elements in $EMM_d(\rho)$, and after a little modification of the proof of Theorem 6.1 in \cite{AguehMartial2011Bitw}, one can show that the barycenter of elements in $EMM_d(\rho)$ is still in $EMM_d(\rho)$. With these properties, we can easily generalize our results to $EMM_d(\rho)$. While \cite{DelonJulie2020Awdi} has a similar generalization, the generalization in our way does not require the marginal consistency. However, our results can not be further generalized to $EMM_d$. In \cite{2024Chen}, it is shown that the barycenter of two elliptically contoured distributions need not to be elliptically contoured. Since the $RW_2$ barycenters are mixtures of barycenters of the elliptically contoured components of the marginal mixtures, their components need not to be elliptically contoured. In other words, the $RW_2$ barycenter of two elements in $EMM_d$ can not belong to $EMM_d$.

\begin{figure}[t]\label{MQ1}
	\includegraphics[width = \textwidth]{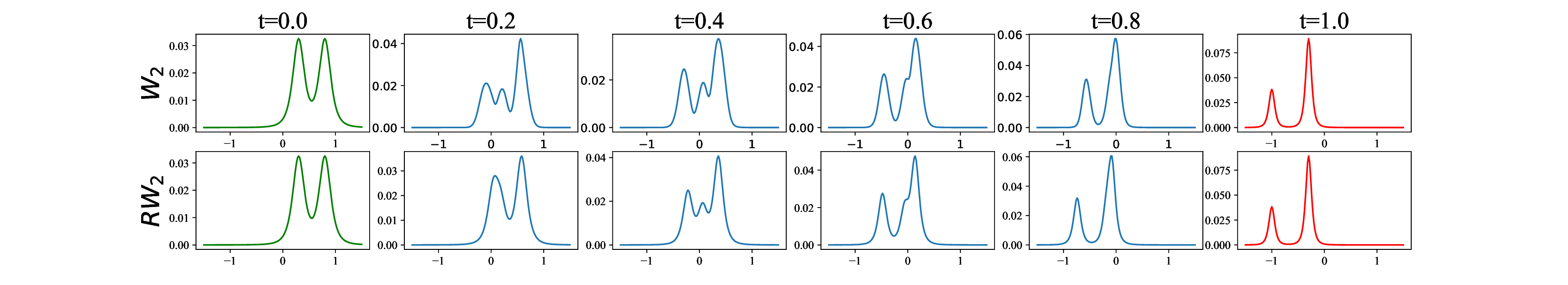}
	\includegraphics[width = \textwidth]{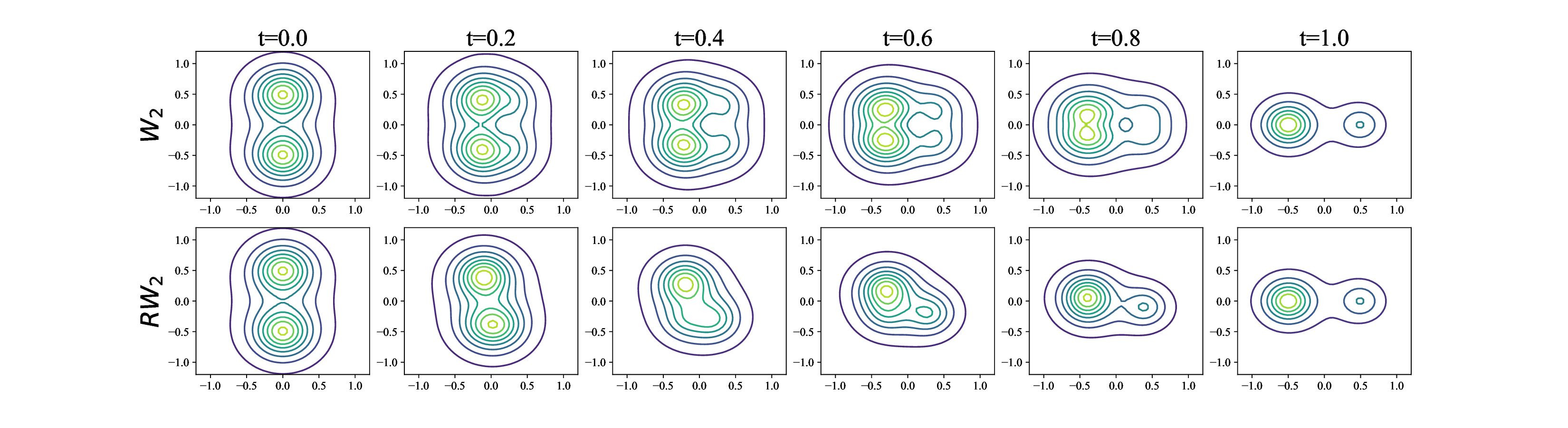}
	\caption{$1$-d and $2$-d examples for $RW_2$ and $W_2$ barycenters of two mixtures of radially contoured distributions with generator $\rho_{IMQ}$. In each case, the first row is the $W_2$ interpolation and the second row is the $RW_2$ interpolation.}
\end{figure}

\section{Some examples for radial mixtures}\label{sec4}

In this section, we compare $RW_2$ barycenters with $W_2$ barycneters using two kinds of generators in $1$-dimensional and $2$-dimensional cases. In each case, we fix the generator $\rho$ and compute the $W_2$ barycenters and $RW_2$ barycenters of given distributions in $RMM_d(\rho)$. The $W_2$ barycenters are all computed using Sinkhorn algorithm with regularization parameter $5 \times 10^{-4}$ in $1$-dimensional case and $5 \times 10^{-3}$ in $2$-dimensional case. The $1$-dimensional case is computed on a grid containing $100$ points, and the $2$-dimensional case is computed on a grid containing $50$ points per dimension.

In terms of the computation cost, it is difficult to find a standard to compare the two methods. The computation cost of $W_2$ barycenters highly depends on the choice of the grid and the regularization parameter, but the computation cost of $RW_2$ barycenters only depends on the number of parameters of radial mixtures. In our setting, the order of magnitude of  an $RW_2$ barycenter is $2$ ms while that of a $1$-dimensional $W_2$ barycenter is $30$ ms and that of a $2$-dimensional $W_2$ barycenter is $3$ s. Therefore, compared to the computation of $RW_2$ barycenters, the computation of $W_2$ is more expensive in several orders of magnitude.

\subsection{Radial mixtures based on inverse multi-quadratic functions}\label{subsec IMQ}

\begin{figure}[t]\label{Compact1}
	\includegraphics[width = \textwidth]{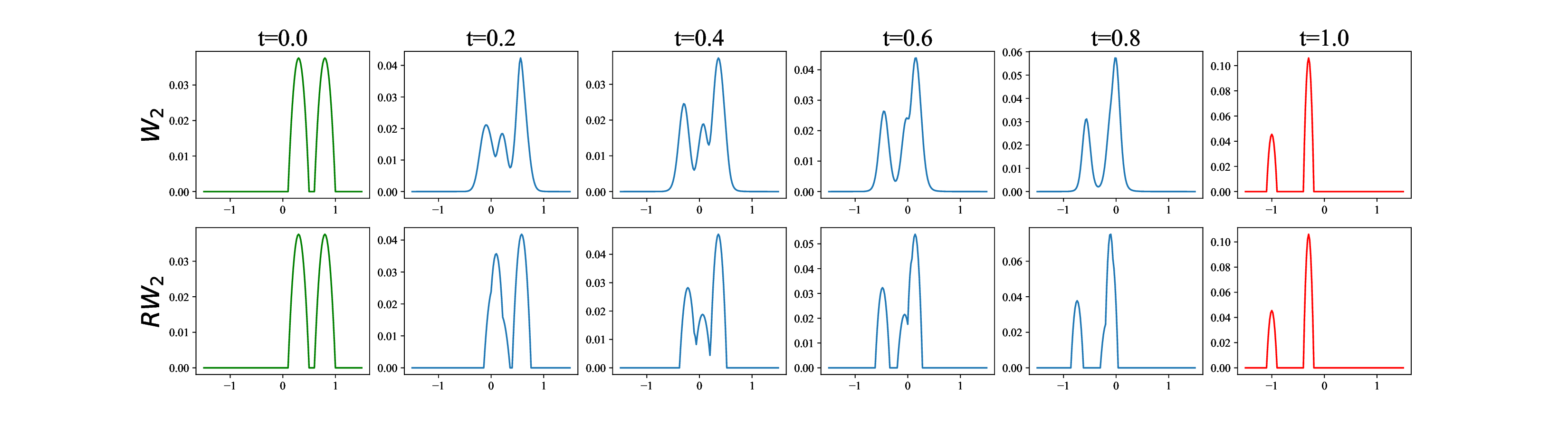}
	\includegraphics[width = \textwidth]{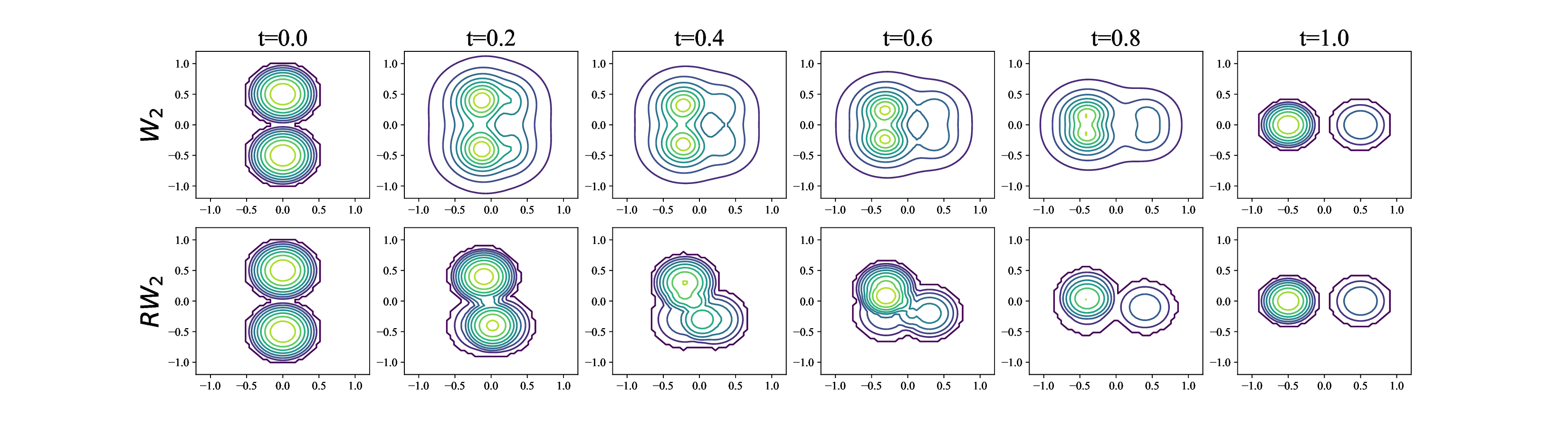}
	\caption{$1$-d and $2$-d examples for $RW_2$ and $W_2$ barycenters of two mixtures of radially contoured distributions with generator $\rho_{C}$. In each case, the first row is the $W_2$ interpolation and the second row is the $RW_2$ interpolation.}
\end{figure}

The multi-quadratic function $(1+x^2)^\alpha$ plays an important role in RBF theory. To ensure such function can be used as a generator, we need $\alpha$ to be negative. That is the inverse multi-quadratic (IMQ) function
\begin{equation}\label{IMQ}
\rho_{IMQ}(x) = (1 + x^2)^{-\beta},\quad \beta > 0.
\end{equation}
To ensure the corresponding distribution is in $\mathcal{P}_2(\R^d)$, it requires that $\beta > d/2$. The radially contoured distribution with generator $\rho_{IMQ}$ is a special $d$-dimensional t-distribution with a scalar covariance matrix. The identifiability of finite mixtures of such distributions has been proved in Example 1 of \cite{HOLZMANNHAJO2006IoFM}. Here we provide $1$-dimensional and $2$-dimensional cases of $W_2$ and $RW_2$ barycenters for mixtures of such radially contoured distributions in \cref{MQ1}. In each case, we set $\beta = 2$ and plot the interpolations induced by Wasserstein distance $W_2$ and $RW_2$.

We observe that the behaviors of the two types of barycenters are quite different. In $1$-d case, the third peak appears later in $RW_2$ barycenters than that appears in $W_2$ barycenters. In $2$-d case, the $W_2$ barycenters are all symmetric about the $x$ axis, while the $RW_2$ barycenters seem like a scaling and rotation of the components.


\subsection{Radial mixtures based on compactly supported functions}

The compactly supported function is another important class of basis functions in RBF theory. The simplest one is the hat function $\rho_{H}(x) = \max\{1-x,0\}$. However, the finite mixture of radially contoured distributions with such generator is not identifiable. It can be found in Example 6 in \cite{HOLZMANNHAJO2006IoFM}. Here we consider the distributions based on such compactly supported function
\begin{equation}
	\rho_{C}(x) = \begin{cases}
		(1 - x^2)^{\frac{1}{\beta}}, \quad & 0 \leq x \leq 1, \\
		0, & x > 1,
	\end{cases} \quad \beta > 0.
\end{equation}
The identifiability of finite mixtures of such distributions can be easily proved in ways similar to Theorem 3 in \cite{HOLZMANNHAJO2006IoFM}. We also give the corresponding examples of $1$-dimensional and $2$-dimensional case in \cref{Compact1}. We set $\beta = 2$ in $1$-dimensional case and $\beta = 3$ in $2$-dimensional case and plot the interpolations induced by Wasserstein distance $W_2$ and our $RW_2$.

The densities of the two marginals are not smooth on the boundaries of their supports due to the rectangular grid. We observe that the densities of $RW_2$ barycenters keep nonsmooth on the boundary, but those of the $W_2$ barycenters do not due to the regularization of Sinkhorn algorithm. Other behaviors are similar to the results presented in the previous section.

\section{Using radial mixtures in practice}\label{sec5}

\subsection{Approximating discrete distributions with radial mixture models}\label{subsec5.1}

An efficient met- hod to approximate discrete distributions by mixture models is using EM algorithm. However, without the explicit formulation as in GMM estimation, there remains an optimization problem in each iteration of the RMM estimation, which makes the computation much slower. To speed up this process, we introduce mini-batch and stochastic technique. The algorithm is illustrated in \cref{mbs EM alg}.

\begin{algorithm}[h]
	\caption{mini-batch stochastic EM for RMM}
	\begin{algorithmic}[1]
		\REQUIRE {\text{dataset} $X = \{x_1, x_2, \cdots, x_N\}$, \text{initial} parameters $\pi^0 = \{\pi_1^0, \pi_2^0, \cdots, \pi_K^0\}$, $m^0 = \{m_1^0, m_2^0, \cdots, m_K^0\}$, $c^0 = \{c_1^0, c_2^0, \cdots, c_K^0\}$, batch size $N_0$.}
		\WHILE{not converged}
		\STATE{\text{Randomly choose} $N_0$ points $\{\bar{x}_1, \bar{x}_2, \cdots, \bar{x}_{N_0}\}$ \text{from} $X$;}
		\FOR{$(i,k) \in \{1, 2, \cdots, N_0\} \times \{1, 2, \cdots, K\}$}
		\STATE $\tau_{i k} = \frac{\pi_k^t \rho(\bar{x}_i; m_k^t, c_k^t)}{\sum_{k=1}^K \pi_k^t \rho(\bar{x}_i; m_k^t, c_k^t)}$;
		\ENDFOR
		\STATE $Q(m, c) = \sum_{i=1}^{N_0} \sum_{k=1}^K \tau_{ik} \log{(\pi_k^t \rho(\bar{x}_n; m_k, c_k))}$;
		\FOR{$k \in \{1, 2, \cdots, K\}$}
		\STATE $\pi_k^{t+1} = \sum_{i=1}^{N_0} \tau_{ik}$
		\ENDFOR
		\STATE $m^{t+1}, c^{t+1} = \operatorname{argmax}_{m,c} Q(m,c)$;
		\STATE $t = t+1$;
		\ENDWHILE
		\RETURN $\pi^t, m^t, c^t$.
	\end{algorithmic}
	\label{mbs EM alg}
\end{algorithm}

Such technique is useful in optimization for large scaling problems \cite{BottouLeon2018OMfL, NemirovskiA2009RSAA}. In the experiments shown in \cref{subsec5.2}, we use \cref{mbs EM alg} with initialized $\pi^0, m^0, c^0$ from K-means algorithm, a batch size $100$ and a maximum number of iterations $1500$. The order of magnitude for getting the parameters is several minutes, which is acceptable in practice. After getting the densities of the radial mixtures, we can use $T_{\rm mean}$ or $T_{\rm rand}$ to replace the position of Monge map in application.

\subsection{Application in color transfer and color averaging}\label{subsec5.2}

We apply our method to the color transfer problem. The color transfer problem requires a source image and a target image. We need to transfer the color style of the source image to that of the target image, so that the output image has the same pattern and geometry as the source image but with the color palette from the target image. Mathematically, an image of size $n_r \times n_c$ can be viewed as a map $u: \Omega \to \R^3$, where $\Omega = \{0, 1, \cdots, n_r - 1\} \times \{0, 1, \cdots, n_c - 1\}$. For each pixel $i \in \Omega$, $u(i) = (r_i, g_i, b_i)$ represents the color of pixel $i$ in RGB space, where the elements $r_i, g_i$ and $b_i$ correspond to the intensities of red, green and blue. In color transfer task, given two images $u_0$ and $u_1$ on grids $\Omega_0$ and $\Omega_1$, define the discrete color distributions of the two images as $\eta_k = \frac{1}{|\Omega_k|} \sum_{i\in\Omega_k} \delta_{u_k(i)}, k = 0, 1$. We need to find a color transfer map $T: \R^3 \to \R^3$ which makes $T_\# \eta_0$ close to $\eta_1$. To utilize our framework, we first approximate the two discrete distributions with two radial mixtures and then compute the transport map. Although our method allows the radial mixtures to come from different families, it will be necessary to solve for the function $T$ in \cref{T} by numerical method in this case, which significantly increases the computation complexity. Therefore, in our experiments, we always use radial mixtures with the same generator. In this experiment, the source image $u_0$ of size $432 \times 576$ and the target image $u_1$ of size $360 \times 480$ are shown in \cref{given figure}.

\begin{figure}[t]\label{error figures}
	\includegraphics[width = \textwidth]{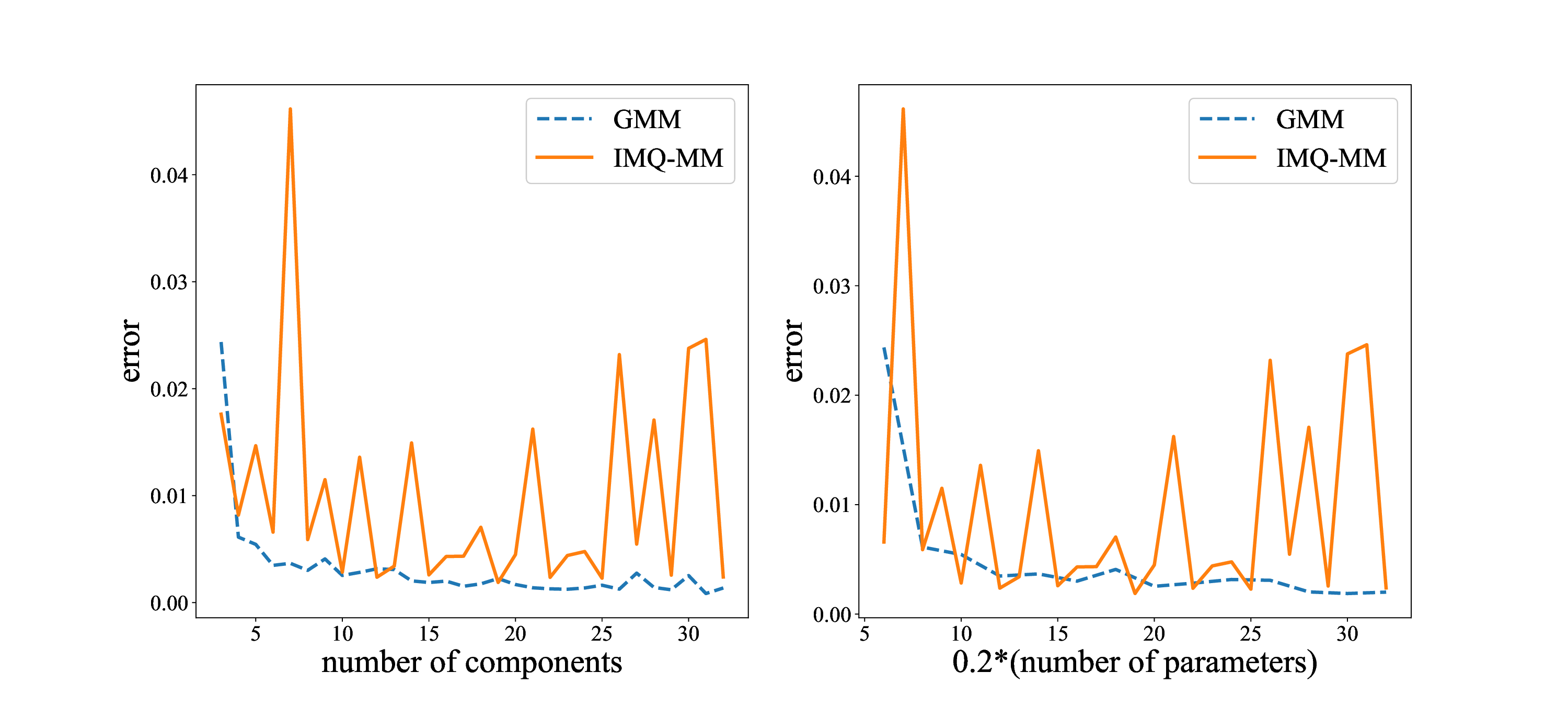}
	\caption{The relationship between error (The square of $W_2$ distance between the distributions of transferred source palette and the target palette) and the number of components(left) and the number of parameters(right). We find that the results of GMM is stable than that of IMQ-MM. When we compare the performance under the same number of components, GMM performs better in most of the time. However, when constrained to the identical number of parameters, the two methods demonstrate comparable results. }
\end{figure}

We will compare our method with that based on GMM proposed in \cite{DelonJulie2020Awdi}. The generator we choose is the inverse multi-quadratic function $\rho_{IMQ}$ given in \cref{subsec IMQ} with $\beta = 3$. The corresponding RMM here we will call IMQ-MM. In both methods, we first use the mixture models to approximate the color distributions of the two images and then compute $T_{\rm mean}$ as the transport map, since $T_{\rm rand}$ will give a noisier result which has been shown in \cite{DelonJulie2020Awdi}. Because $T_{\rm mean}$ is not a measure preserving map as discussed in \cref{subsec33}, the first thing we care is how far it is between $T_{\rm mean}$ and a transport map. We use the following value as a principle:

\begin{figure}[t]\label{given figure}
	\includegraphics[width=0.49\textwidth]{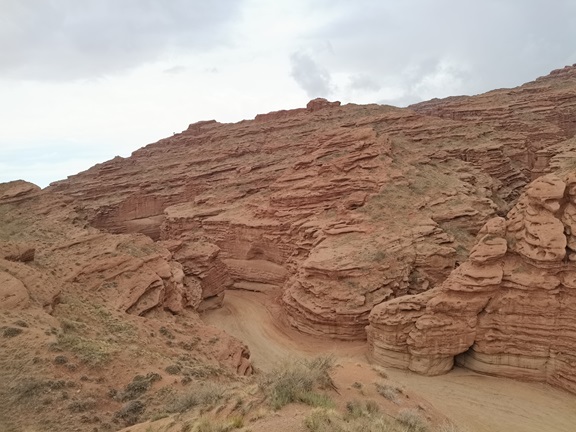}
	\includegraphics[width=0.49\textwidth]{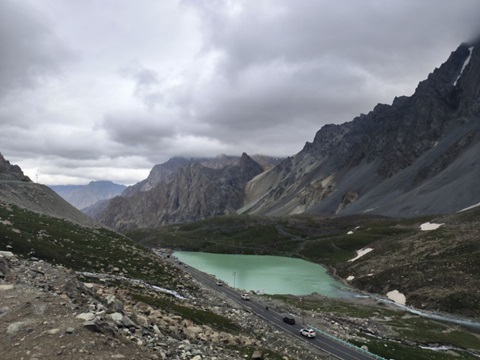}
	
	\caption{Left, the source image $u_0$. Right, the target image $u_1$.}
\end{figure}

\begin{equation}\label{error}
	\text{error} = W_2((T_{\rm mean})_\# \eta_0, \eta_1)^2.
\end{equation}
If $T_{\rm mean}$ is a transport map, we have $(T_{\rm mean})_\# \eta_0 = \eta_1$, and the error will be $0$. Intuitively, while increasing the number of components, approximation will be better. We plot error vs the number of components on the left of \cref{error figures}. On the other hand, a typical advantage of RMM over GMM is that when we fix the number of components, $n$, GMM needs to fit one weight vector, $n$ mean vectors, and $n$ covariance matrices, which correspond to $(1+1.5d+0.5d^2)n$ parameters, where $d$ is the dimension of the underlying space, but RMM only needs to fit $(d+2)n$ parameters. The color transfer task is a $3$-dimensional problem. Thus the number of parameters in GMM is $10n$, whereas that in RMM is $5n$. We also plot the error vs the number of parameters on the right of \cref{error figures}. Here the color distributions of the images are all approximated by mixture models with the same number of components. Note that the OT problem \cref{error} is too large to compute. Thus, in this experiment, we randomly choose $m$ pixels in each image as $\eta_0^{m}$ and $\eta_1^{m}$ and compute the $W_2$ distance between $(T_{\rm mean})_\# \eta_0^m$ and $\eta_1^m$. Here we set $m = 5000$ and repeat the process $50$ times. Each error value shown in \cref{error figures} is the average of the results computed $10$ times. The theoretical fundamental of this operation is based on the results of \cite{DudleyR.M.1969TSoM}, which reads
\begin{equation*}
	E\left(\left|W_2(\hat{\alpha}_n, \hat{\beta}_n) - W_2(\alpha, \beta)\right|\right) = O\left(n^{-\frac{1}{d}}\right),
\end{equation*}
where $\hat{\alpha}_n$ and $\hat{\beta}_n$ are discrete distributions of $n$ samples from compactly supported probability distributions $\alpha$ and $\beta$.

\begin{figure}[t]\label{transferred figure}
	\includegraphics[width=0.32\textwidth]{fig3.jpg}
	\includegraphics[width=0.33\textwidth]{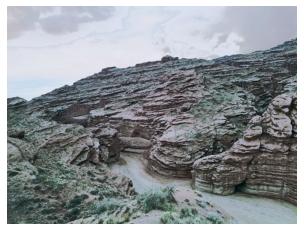} 
	\includegraphics[width=0.33\textwidth]{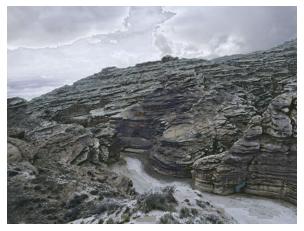}
	\includegraphics[width=\textwidth]{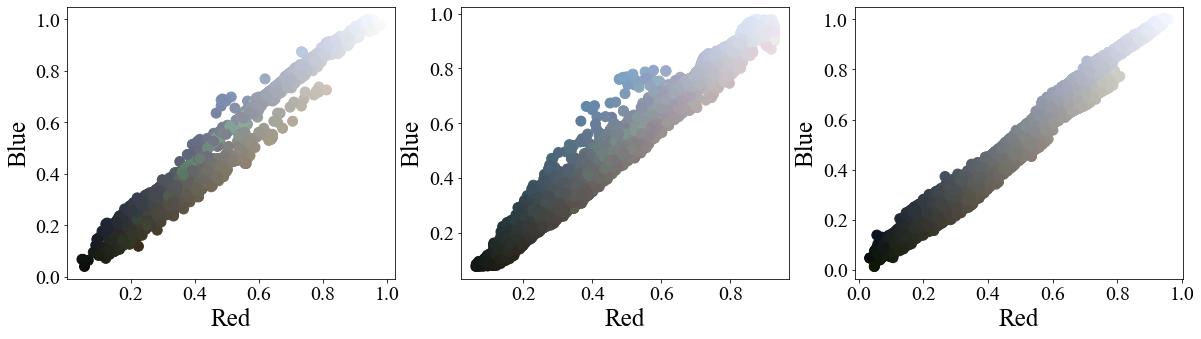}
	\includegraphics[width=\textwidth]{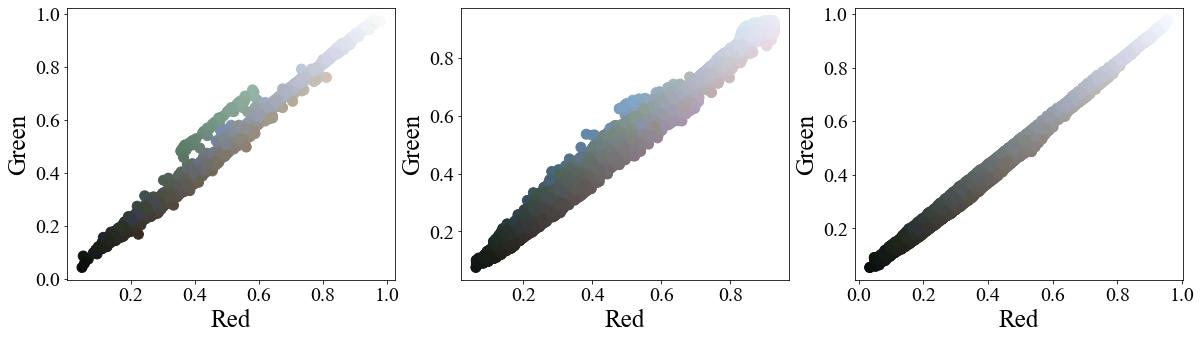}
	\caption{First row: left, the target image. Middle, the transferred image by IMQ-MM. Right, the transferred image by GMM. Second row, the projections of color palettes of images to blue and red space. Left, the target image. Middle, the top left image. Right, the top right image. Third row, the projections of color palettes of images to red and green space. The order is the same.}
		
\end{figure}

\cref{error figures} shows that the error of GMM is more stable than that of IMQ-MM. It can be attributed to the elliptical shape of GMM components, which enables superior adaptability to practical data distribution, whereas the radial configuration of the components of IMQ-MM renders its fitting performance more sensitive to the number of components. In the left figure, better performance of GMM over IMQ-MM in most scenarios is partially attributable to its utilization of twice as many parameters as in IMQ-MM. When constrained to an identical parameter count, the two models show comparable performance, as is revealed in the right figure. Furthermore, our results indicate that the increasing of the number of components does not significantly reduce the error, but increases the computation complexity. After a trade-off between the error performance and computational efficiency, we choose $10$ components for GMM and $15$ components for IMQ-MM. The results are shown in the first row of \cref{transferred figure}.

We find that the image transferred by IMQ-MM is lighter than by GMM. Both images successfully keep the pattern and geometry of the source image. We further consider the differences of the color distributions between the target image and the two transferred images. We plot the projections of these color distributions on RGB space to B\&R space and G\&R space in the last two rows in \cref{transferred figure}. It is clear that GMM gives a ``thinner" distribution than the target distribution. Especially in the G\&R space, some points of the target distribution are lost. However, the distribution given by IMQ-MM is ``fatter" than the target distribution. Both results have some differences from the target distribution, but IMQ-MM seems better in this experiment.

\begin{figure}[t]\label{averaging figures}
	\centering
	\includegraphics[width=0.24\textwidth]{fig3.jpg}
	\includegraphics[width=0.24\textwidth]{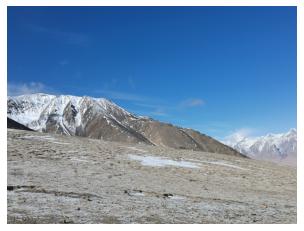}
	\rotatebox{90}{\includegraphics[width=0.24\textwidth]{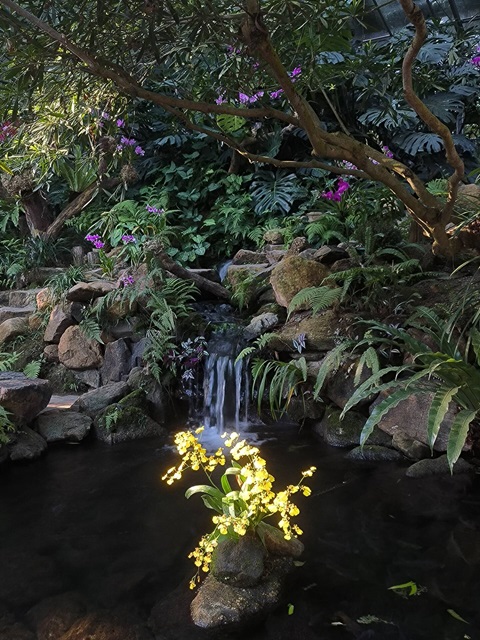}} \\
	\includegraphics[width=0.24\textwidth]{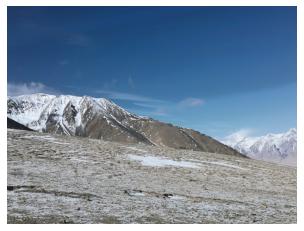}
	\includegraphics[width=0.24\textwidth]{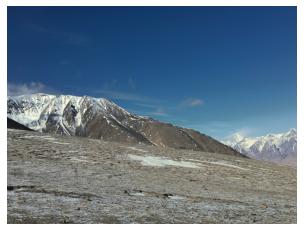} \\
	\includegraphics[width=0.24\textwidth]{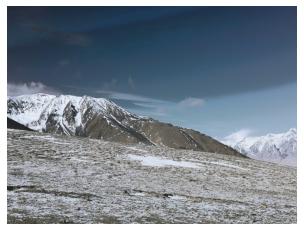}
	\includegraphics[width=0.24\textwidth]{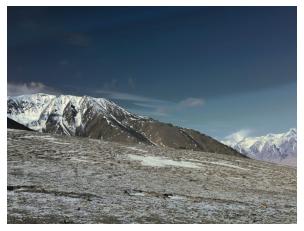}
	\includegraphics[width=0.24\textwidth]{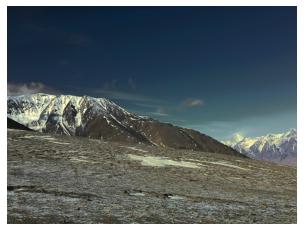} \\
	\includegraphics[width=0.24\textwidth]{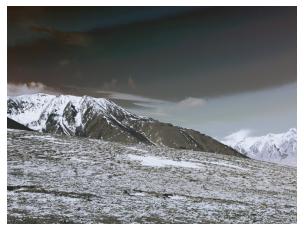}
	\includegraphics[width=0.24\textwidth]{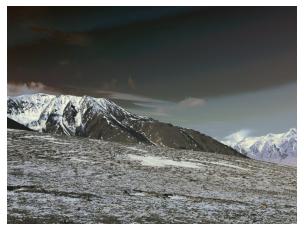}
	\includegraphics[width=0.24\textwidth]{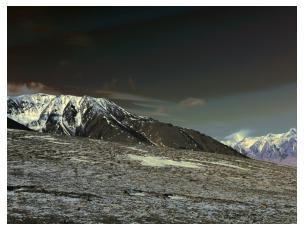}
	\includegraphics[width=0.24\textwidth]{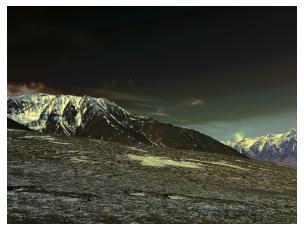}
	\caption{In this experiment, the top middle image is modified in such way that its color distribution goes through the IMQ-$W_2$ barycenters of the color distributions of the top $3$ images. Each source is approximated by an IMQ-mixture of $14$ components. The weights used in the barycenters are the barycenter coordinates with respect to the top $3$ images of the top $3$ images in the triangle. The weight of the bottom left one is $1$ of the top left image and that of the bottom right one is $1$ of the top right image.}
\end{figure}

In terms of the cost of computation, our method based on IMQ function is slower than GMM. The order of magnitude for outputting the top middle image in \cref{transferred figure} is $5$ minutes and that for outputting the top right image is $50$ seconds. The main cost is spent in the EM algorithm in this experiment. When we use EM algorithm to update the parameters, there exists an explicit formulation in each step of GMM. However, each updating step in IMQ-MM requires us to solve an optimization problem.

We end this section with a color averaging experiment showed in \cref{averaging figures}. Three images are given at the top. The color palettes of the images are represented by RMMs with $14$ components each. We show the images transferring the color palette of the top middle image to the $RW_2$ barycenters of the $3$ color palettes. The barycenter weight of the bottom left corner is $1$ of the top left image and that of the bottom right one is $1$ of the top right image. The weights of other images are the corresponding barycenter coordinates. For example, the weight of the second image in the third line is $(1/3, 1/3, 1/3)$.

\section{Conclusions and remarks}\label{sec6}

In this paper, we define a relaxed Wasserstein distance $RW_2$ on the set of all identifiable RMMs by viewing the RMMs as discrete distributions in the space of radially contoured distributions. We show that $RMM_d$ equipped with $RW_2$ is a geodesic space, and the computation of transport plan and barycenter in this space is quite easier than in Wasserstein space. Compared to $GW_2$ defined in \cite{DelonJulie2020Awdi}, which can only be generalized to $EMM_d(\rho)$ whose components satisfy marginal consistency, our $RW_2$ can be generalized to all the $EMM_d(\rho)$. In application, we use EM algorithm to estimate the mixtures from discrete data. Since the M-step of the EM algorithm for RMMs requires us to solve an optimization problem while that for GMMs has an explicit formulation, the EM algorithm for RMMs is much slower than that for GMMs. In \cref{subsec5.1}, we introduce the mini-bath stochastic EM \cref{mbs EM alg} to speed up estimating RMMs. In an experiment of color transfer task, we have compared our method with that proposed in \cite{DelonJulie2020Awdi}. In this task, two methods provide comparable numerical results, and the color distribution of the output image is more desirable by our method.

This work can be extended in several ways in the future. On the one hand, for example, a relaxed Gromov-Wasserstein distance analogous to $RW_2$ can be developed. On the other hand, if we allow the coefficients to be negative, a continuous way to transfer a function to another can be established by fitting each function with RBFs and continuously changing the parameters.

%

\section*{Acknowledgments}
We would like to thank Songyan Luo for the helpful discussion.

\bibliographystyle{siamplain}
\bibliography{references}
\end{document}